\documentclass[10pt,a4paper,leqno]{article}
\usepackage[TS1,T1]{fontenc}
\usepackage{textcomp,lmodern}
\input{ts1lmr.fd}
\usepackage{amsmath,amssymb,amsthm,mathtools}
\usepackage[margin=3cm]{geometry}
\usepackage{titlesec,fancyhdr}
\usepackage{color,xcolor}
\usepackage[hidelinks]{hyperref}
\titleformat{\subsection}{\normalfont\itshape}{\thesubsection.\enspace}{1pt}{}
\numberwithin{equation}{section}
\newtheorem{theorem}{Theorem}[section]
\newtheorem{lemma}[theorem]{Lemma}
\newtheorem{proposition}[theorem]{Proposition}
\theoremstyle{remark}
\newtheorem{remark}[theorem]{Remark}
\newtheorem*{remark*}{Remark}
\newcommand{\R}{\mathbb R}
\newcommand{\dd}{\,\mathrm d}
\newcommand{\dt}{\frac{\mathrm d}{\mathrm dt}}
\newcommand{\proofstep}[1]{\par\smallskip\noindent\textbf{#1}\par\nobreak}
\allowdisplaybreaks[2]
\title{A Sharp Criterion for Finite-Time Blow-Up and Global Existence
in the Degasperis--Procesi Equation}
\author{Zhaonan Luo$^{1}$\,\thanks{Email: \texttt{luozhn7@mail.sysu.edu.cn}.},\hspace{0.7em}
Weikui Ye$^{2}$\,\thanks{Email: \texttt{yeweikuisz@szu.edu.cn}.}\hspace{0.7em}and
Zhaoyang Yin$^{1}$\,\thanks{Email: \texttt{mcsyzy@mail.sysu.edu.cn}.}\\[0.7em]
{\small $^{1}$\,School of Science, Shenzhen Campus of Sun Yat-sen University,}\\[-0.2em]
{\small Shenzhen 518107, China}\\[0.2em]
{\small $^{2}$\,School of Mathematical Sciences, Shenzhen University,}\\[-0.2em]
{\small Shenzhen 518060, China}}
\date{}
\begin{document}
\maketitle
\hrule

\begin{abstract}
The Camassa--Holm and Degasperis--Procesi equations are two prototypical
models for nonlinear shallow-water waves. Determining whether a given
smooth initial profile gives rise to a global strong solution or to
finite-time blow-up is a central problem in the analysis of these
equations, even for other nonlinear shallow-water waves. For the Camassa--Holm equation, Constantin and Escher's
wave-breaking analysis \cite{ConstantinEscher1998}, McKean's classification
in terms of the initial momentum \cite{McKean2004}, and Brandolese's
single-point blow-up criterion \cite{Brandolese2014}, among others, have made fundamental
contributions to the theory of blow-up and global existence. For the
Degasperis--Procesi equation, prior results give criteria for particular
classes of initial data (see, e.g., \cite{LiuYin2006}), leaving the general
classification open. In this paper, we establish the following necessary
and sufficient condition for finite-time blow-up of solutions with smooth
initial data:
$$T^*<\infty
\quad\Longleftrightarrow\quad
\exists\,a<b:\quad m_0(a)>0>m_0(b).$$
This resolves the problem of classifying finite-time blow-up and global
existence for the Degasperis--Procesi equation.
\end{abstract}


\tableofcontents

\section{Introduction}

\subsection{The equation and previous results}

We consider the Cauchy problem for the Degasperis--Procesi (DP) equation
\begin{equation}\label{eq:dp-local}
\begin{cases}
u_t-u_{txx}+4uu_x=3u_xu_{xx}+uu_{xxx},&t>0,\ x\in\R,\\
u(0,x)=u_0(x),&x\in\R.
\end{cases}
\end{equation}
Its nonlocal and momentum formulations are, respectively,
\begin{equation}\label{eq:dp}
u_t+uu_x=-\Pi_x,\qquad
\Pi=\frac32p*u^2,\qquad p(x)=\frac12e^{-|x|},
\end{equation}
and
\begin{equation}\label{eq:momentum}
m=u-u_{xx},\qquad m_t+um_x+3u_xm=0.
\end{equation}
Here $*$ denotes convolution in the spatial variable. The local equation
and the momentum equation are understood in the sense of distributions
whenever the regularity does not permit a pointwise interpretation.

For $u_0\in H^s(\R)$, $s>3/2$, let $u$ be the maximal strong solution,
with maximal existence time $T^*$:
\begin{equation}\label{eq:solution}
u\in C([0,T^*);H^s(\R))
\cap C^1([0,T^*);H^{s-1}(\R)),\qquad T^*\in(0,\infty].
\end{equation}
Degasperis and Procesi \cite{DegasperisProcesi1999} identified the DP
equation through the method of asymptotic integrability.
Degasperis, Holm, and Hone \cite{DegasperisHolmHone2002} subsequently
constructed a Lax pair, proposed a bi-Hamiltonian structure, and
derived conservation laws and peakon solutions.
Constantin and Lannes \cite{ConstantinLannes2009} justified the relevance
of the Camassa--Holm and DP equations to the propagation of shallow
water waves over a flat bed, in a regime that retains stronger nonlinear
effects than the classical Korteweg--de Vries approximation.
The combination of integrability and blow-up makes the DP
equation a useful model for studying the competition between nonlinear
steepening and nonlocal interactions.

Firstly, we introduce the well-known  Camassa-Holm (CH) equation
\begin{align}\tag{CH}
u_{t}-u_{txx}+3uu_{x}=2u_{x}u_{xx}+uu_{xxx}.
\end{align}
The CH equation has been studied extensively in the last three decades because of its many remarkable distinctive properties. One of the particular features is that it has single peakon solutions and multipeakon solutions. The local well-posedness and ill-posedness for the CH equation in Sobolev and Besov spaces were established in \cite{d2,lo}. Furthermore, global strong solutions \cite{ce4} and global weak solutions \cite{chk} have been established. Another particular feature is the blow-up phenomena (the solutions remain bounded but their slopes become unbounded in finite time). In 1998, Constantin and Escher
\cite{ConstantinEscher1998} developed a rigorous method for tracking
the minimum slope of a solution and showed that a sufficiently steep
negative initial slope leads to blow-up for the Camassa--Holm
equation. McKean \cite{McKean1998,McKean2004} gave a landmark
classification on the real line: for $u_0\in H^3(\mathbb{R})$ with
integrable initial momentum $m_0=u_0-u_0''$, blow-up occurs in
finite time if and only if $m_0$ is positive at some point to the left
of a point where it is negative. Furthermore, Brandolese \cite{Brandolese2014}
established the elegant single-point criterion
$u_0'(x_0)<-|u_0(x_0)|$ for finite-time blow-up in the CH equation. His Lyapunov-function
argument unified several earlier sufficient conditions and extended
to shallow-water-type equations beyond the integrable setting.

Now, we present the Cauchy problem for the DP equation in detail.
Yin \cite{Yin2003} established local well-posedness in $H^s(\R)$,
$s>3/2$, persistence of regularity, and a continuation criterion,
while Gui and Liu \cite{GuiLiu2011DP} developed the Besov theory.
Further results concern nonuniform dependence on initial data
\cite{HimonasHolliman2011DP}, norm inflation at critical regularity
\cite{GuoLiuMolinetYin2019}.
The DP equation also admits peakons $u(t,x)=ce^{-|x-ct|}$.
Using a Lyapunov functional, Lin and Liu \cite{LinLiu2009} proved
orbital stability of the positive-speed peakons for perturbations
whose initial momentum is a nonnegative finite Radon measure.

Global existence and finite-time breakdown are closely connected to
the momentum $m_0=u_0-u_0''$.
Yin \cite{Yin2004} constructed global strong and weak solutions on the
real line under suitable sign assumptions on $m_0$.
Early discussions of global existence and blow-up also include
Zhou \cite{Zhou2004DP}.
Liu and Yin \cite[Theorems~3.1 and~3.2]{LiuYin2006} proved global
existence for $u_0\in H^s(\R)$, $s>3/2$, when $m_0$ has one sign,
or when there exists $x_*\in\R$ such that
\[
m_0(x)\le0\quad(x\le x_*),\qquad
m_0(x)\ge0\quad(x\ge x_*).
\]
The corresponding slope bounds are $|u_x|\le|u|$ in the same-sign
case and $u_x\ge-|u|$ in the second case.
For the opposite sign arrangement, Liu and Yin
\cite[Theorem~4.2]{LiuYin2006} proved finite-time blow-up under the
conditions
\[
m_0(x)\ge0\quad(x\le x_*),\qquad
m_0(x)\le0\quad(x\ge x_*),\qquad m_0\text{ changes sign}.
\]
Their later results \cite[Theorems~3.1 and~3.2]{LiuYin2007} allow
odd momentum profiles with additional sign changes: for
$u_0\in H^s(\R)$, $s>3/2$, blow-up occurs if $m_0$ is odd and,
for some $a>0$, $m_0(-a)=0$ and either
\[
\begin{cases}
m_0(x)>0,&x<-a,\\
m_0(x)<0,&-a<x<0,
\end{cases}
\qquad\text{or}\qquad
\begin{cases}
m_0(x)\le0,&x<-a,\\
m_0(x)>0,&-a<x<0.
\end{cases}
\]
Escher, Liu, and Yin \cite[Theorem~3.3]{EscherLiuYin2006DP}
showed that blow-up occurs only at the origin when
\[
u_0\in H^s(\R),\quad s>\frac32,\quad
u_0\not\equiv0,\quad u_0(-x)=-u_0(x),\quad
m_0(x)\ge0\quad(x<0).
\]
They also studied the blow-up rate and constructed global weak solutions.

A complementary approach seeks conditions involving the velocity and
its derivative at a single point. Without a symmetry assumption, Liu and Yin
\cite[Theorem~4.1]{LiuYin2006} obtained blow-up if, for some
$\varepsilon>0$ and $x_0\in\R$,
\[
u_0'(x_0)<-\frac{(1+\varepsilon)\sqrt6}{4}
\left[\|u_0\|_\infty+
\sqrt{\|u_0\|_\infty^2+
2\sqrt6\,\|u_0\|_2^2\log\left(1+\frac2\varepsilon\right)}\right].
\]
Chen, Guo, Liu, and Qu \cite[Theorem~2.1]{ChenGuoLiuQu2016} obtained
a sufficient condition involving the conserved quadratic energy:
\[
u_0'(x_0)<-\sqrt{\frac32}|u_0(x_0)|
-\sqrt{3\left(\sqrt{\frac32}-1\right)E_0},
\qquad x_0\in\R,
\]
where
\[
E_0=\int_\R u_0(1-\partial_x^2)(4-\partial_x^2)^{-1}u_0\dd x.
\]
In the same Sobolev range, Zheng \cite[Theorem~3.3]{Zheng2025}
proved a pointwise criterion in the smoothed velocity
\[
v_0=(4-\partial_x^2)^{-1}u_0:
\qquad 3v_0'(x_0)<-\left|2v_0(x_0)+v_0''(x_0)\right|
\quad\text{for some }x_0\in\R.
\]
Specializing the result of Escher and Yin
\cite[Theorem~5.2]{EscherYin2008b} to $b=3$ with zero linear
dispersion gives finite-time blow-up for the DP equation when
\[
u_0\in H^s(\R),\quad s>\frac32,\quad
u_0\not\equiv0,\quad u_0(-x)=-u_0(x),\quad u_0'(0)\le0.
\]
For the periodic DP equation, global strong solutions and
finite-time blow-up have also been studied; see
\cite{Yin2003PeriodicDP,Yin2004PeriodicDP,EscherLiuYin2007PeriodicDP,BrandoleseDP2016}.


\subsection{Main results}
The works cited above establish criteria for finite-time blow-up and
global existence for several classes of initial data for the
Degasperis--Procesi (DP) equation. This leads to a natural question:

\textbf{
Given an arbitrary smooth initial datum, can one determine from that
datum alone whether the solution exists globally or
blows up in finite time?}

 This question belongs to the broader study
of global regularity and singularity formation in nonlinear partial
differential equations. A prominent example is the Millennium Prize
Problem on the existence and smoothness of solutions to the
three-dimensional incompressible Navier--Stokes equations. In this
paper, we answer the question for the DP equation by establishing
a necessary and sufficient condition on the initial data for global
existence of the corresponding solution.

\begin{theorem}[Momentum characterization of blow-up]\label{thm:momentum}
Let $u_0\in H^s(\R)$, $s>5/2$, and set $m_0=u_0-u_0''$.
Then the maximal strong solution of \eqref{eq:dp} satisfies
\begin{equation}\label{eq:momentum-iff}
T^*<\infty
\quad\Longleftrightarrow\quad
\exists\,a<b:\quad m_0(a)>0>m_0(b).
\end{equation}
Equivalently, the solution is global if and only if $m_0\ge0$ on $\R$,
{or} $m_0\le0$ on $\R$, or there exists $x_*\in\R$ such that
\begin{equation}\label{eq:global-sign}
m_0(x)\le0\quad(x\le x_*),\qquad
m_0(x)\ge0\quad(x\ge x_*).
\end{equation}
If $T^*<\infty$, blow-up occurs in the sense that
\begin{equation}\label{eq:breaking}
\sup_{0\le t<T^*}\|u(t)\|_\infty<\infty,
\qquad \lim_{t\uparrow T^*}\inf_{x\in\R}u_x(t,x)=-\infty.
\end{equation}
\end{theorem}

For initial data of lower regularity, we give a pointwise criterion for blow-up.
\begin{theorem}[A sharp pointwise criterion for blow-up]\label{thm:single}
Let $u_0\in H^s(\R)$, $s>3/2$. If there exists $x_0\in\R$ such that
\begin{equation}\label{eq:single-initial}
u_0'(x_0)<-C|u_0(x_0)|,\quad C=1,
\end{equation}
then the solution of \eqref{eq:dp} blows up in finite time and the constant $C$ here is sharp.
\end{theorem}

\begin{theorem}[A finite-time characterization of blow-up]\label{thm:entry}
Let $u_0\in H^s(\R)$, $s>3/2$, and let $u$ be the maximal strong
solution of \eqref{eq:dp}. Then
\begin{equation}\label{eq:entry-iff}
T^*<\infty
\quad\Longleftrightarrow\quad
\exists\,t_1\in[0,T^*),\ x_1\in\R:
\quad u_x(t_1,x_1)<-|u(t_1,x_1)|.
\end{equation}
Whenever these equivalent conditions hold, \eqref{eq:breaking} holds.
\end{theorem}

\begin{remark}[Relation between the momentum and pointwise conditions]\label{rem:criteria-relation}
Let $u_0\in H^s(\R)$ with $s>5/2$, and suppose that
\eqref{eq:single-initial} holds at $x_0$. By \eqref{eq:momentum-convolution},
\[
\begin{aligned}
\int_{-\infty}^{x_0}e^{y-x_0}m_0(y)\dd y
&=u_0(x_0)-u_0'(x_0)>0,\\
\int_{x_0}^{\infty}e^{x_0-y}m_0(y)\dd y
&=u_0(x_0)+u_0'(x_0)<0.
\end{aligned}
\]
Since the weights are positive, there exist $a<x_0<b$ such that
$m_0(a)>0>m_0(b)$. Thus the initial pointwise condition implies the
momentum sign condition at this regularity. By
Theorems~\ref{thm:momentum} and~\ref{thm:entry}, the momentum sign
condition is equivalent to the pointwise condition holding at
some time $t_1\in[0,T^*)$. This implies that $\exists\,t_1\in[0,T^*),\ x_1\in\R:
\quad u_x(t_1,x_1)<-|u(t_1,x_1)|$ is a necessary and sufficient condition for finite-time blow-up of the DP solution. However, when only the initial data are considered, $\exists\,x_1\in\R:
\quad u_x(0,x_1)<-|u(0,x_1)|$ remains a sufficient condition for finite-time blow-up of the DP solution.
\end{remark}

\begin{remark}[An $H^2$ example]
To illustrate the sharp pointwise criterion of Theorem~\ref{thm:single} at
$H^2$ regularity, consider the piecewise rational datum
\[
u_0(x)=
\begin{cases}
\dfrac{1-2x}{1+x^2},&x<0,\\[1mm]
\dfrac{1-2x}{1+2x^2},&x\ge0.
\end{cases}
\]
The two pieces match in value and first derivative at zero, with
$u_0(0)=1$ and $u_0'(0)=-2$. Thus the weak second derivative agrees
with the classical second derivatives on the two half-lines. Since
$u_0^{(j)}(x)=O(|x|^{-1-j})$ as $|x|\to\infty$ for $j=0,1,2$,
we have $u_0\in H^2(\R)$.
The one-sided limits $u_0''(0^-)=-2$ and $u_0''(0^+)=-4$ are distinct.
For sufficiently small $\varepsilon>0$, the square of the
$H^{1/2}$ seminorm of $u_0''$ is therefore bounded below by a
positive multiple of
\[
\int_0^\varepsilon\int_0^\varepsilon\frac{\dd x\,\dd y}{(x+y)^2}
=+\infty.
\]
Hence $u_0\notin H^{5/2}(\R)$ and lies outside the regularity range of
Theorem~\ref{thm:momentum}. Since $u_0'(0)=-2<-|u_0(0)|=-1$,
Theorem~\ref{thm:single} gives $T^*<\infty$, and
Theorem~\ref{thm:entry} applies with $t_1=x_1=0$. It appears to be an interesting problem to establish a dichotomy between blow-up and global existence in terms of the initial data alone at low regularity $(u_0\in H^s(\R),~s>3/2)$.
\end{remark}

\subsection{Motivation and strategy}

The proof of Theorem~\ref{thm:momentum} uses arguments from
the proofs of Theorems~\ref{thm:single} and~\ref{thm:entry}. We obtain the
pointwise criterion and the finite-time characterization by studying
a spatial maximum. Set
\begin{equation}\label{eq:intro-AB}
A(t,x)=u(t,x)-u_x(t,x),\qquad
B(t,x)=-u(t,x)-u_x(t,x).
\end{equation}
Then $u_x<-|u|$ is equivalent to $A>0$ and $B>0$.
Let $q(t,a)$ be the characteristic starting at $a$,
\begin{equation}\label{eq:intro-flow}
q_t(t,a)=u(t,q(t,a)),\qquad q(0,a)=a.
\end{equation}
For $\rho>0$ and bounded continuous $f$, define
\begin{equation}\label{eq:intro-K}
\begin{aligned}
K_{\rho,-}f(x)&=\int_{-\infty}^x e^{\rho(y-x)}f(y)\dd y,\\
K_{\rho,+}f(x)&=\int_x^\infty e^{\rho(x-y)}f(y)\dd y,
\end{aligned}
\end{equation}
and write $K_\pm=K_{1,\pm}$.
 
For the Camassa--Holm equation, the slope along a characteristic satisfies
   $$\frac{\mathrm d}{\mathrm dt}u_x(t,q)\leq\frac12(u^2-u^2_x)(t,q).$$
McKean \cite[Section~2, formulas (3) and (4)]{McKean1998}
established the following weighted identities:
\[
\begin{aligned}
\frac{\mathrm d}{\mathrm dt}\bigl[e^qA(t,q)\bigr]
&=\frac12e^q\bigl[A^2+K_-(A^2)\bigr]\ge0,\\
\frac{\mathrm d}{\mathrm dt}\bigl[e^{-q}B(t,q)\bigr]
&=\frac12e^{-q}\bigl[B^2+K_+(B^2)\bigr]\ge0.
\end{aligned}
\]
If $A(0,a)>0$ and $B(0,a)>0$, these identities imply that
$(AB)(t,q(t,a))$ is nondecreasing. Together with the preceding
differential inequality, this yields the pointwise blow-up criterion.

However, for the DP equation,  we only deduce that
\begin{equation}\label{eq:intro-DP-raw-AB-weighted}
\begin{aligned}
\frac{\mathrm d}{\mathrm dt}\bigl[e^{2q}A(t,q)\bigr]
&=e^{2q}\left[A^2+\frac12K_-(A^2)-\frac12K_-(AB)\right],\\
\frac{\mathrm d}{\mathrm dt}\bigl[e^{-2q}B(t,q)\bigr]
&=e^{-2q}\left[B^2+\frac12K_+(B^2)-\frac12K_+(AB)\right].
\end{aligned}
\end{equation}
The additional terms involving $K_\pm(AB)$ have no prescribed sign.
These identities do not ensure that $A$ and $B$ remain positive
along a fixed characteristic, and therefore do not yield the
same monotonicity argument for the DP equation.

To overcome this difficulty, we consider \textbf{a new quantity $J(t)$}
to track the pointwise condition for blow-up, defined by
\begin{equation}\label{eq:intro-J}
J(t)=\sup_{x\in\R}\bigl(\min\{A(t,x),B(t,x)\}\bigr)_+
    =\sup_{x\in\R}\bigl(-u_x(t,x)-|u(t,x)|\bigr)_+.
\end{equation}
Here $r_+=\max\{r,0\}$, so $J(t)>0$ precisely when $u_x<-|u|$
somewhere at time $t$. Since $u_x(t,x)\to0$ as $|x|\to\infty$,
the definition gives
\begin{equation}\label{eq:J-slope}
J(t)\le-\inf_{x\in\R}u_x(t,x)\le J(t)+\|u(t)\|_\infty,
\end{equation}
which implies that $J(t)\rightarrow+\infty ~\Longrightarrow~
\inf_{x\in\R}u_x(t,x)\rightarrow-\infty$. 

At a differentiability time $t$ with $J(t)>0$, choose $a$ such that
\[
J(t)=-u_x(t,q(t,a))-|u(t,q(t,a))|.
\]
In the case $u(t,q(t,a))>0$, we obtain
\begin{align}
J'(t)
 \ge \big(u^2_x-\frac32u^2 +\frac32K_+u^2\big)(t,q(t,a))\ge u^2_x(t,q(t,a))-\frac32u^2 (t,q(t,a)).
\end{align}
A natural approach is to drop the nonnegative term
$\frac{3}{2}K_+(u^2)(t,q(t,a))$ and use some conservation laws to control
the negative term $-\frac{3}{2}u^2(t,q(t,a))$. Discarding this favorable
contribution weakens the resulting blow-up criterion and leaves out a
substantial class of initial data that lead to finite-time blow-up.
To retain this contribution, we derive a pointwise lower bound for
$\frac{3}{2}K_+(u^2)(t,q(t,a))$ that offsets part of the negative term
$-\frac{3}{2}u^2(t,q(t,a))$, yielding the differential inequality below:
\[
J'(t)\ge J(t)^2\qquad\text{for a.e. }t\in(0,T^*).
\]
This proves the pointwise blow-up criterion in
Theorem~\ref{thm:single}.

Next, to prove the wave-breaking implication of Theorem~\ref{thm:momentum}, assume
$m_0(a)>0>m_0(b)$ for some $a<b$, and suppose that the solution is
global. We first use the one-sided smoothed quantities and their representations:
\begin{gather}
\begin{aligned}
\mathcal A(t,x)&=K_{2,-}A(t,x)
=\int_{-\infty}^x e^{2(y-x)}(u-u_x)(t,y)\dd y,\\
\mathcal B(t,x)&=K_{2,+}B(t,x)
=\int_x^\infty e^{2(x-y)}(-u-u_x)(t,y)\dd y,
\end{aligned}\label{eq:intro-filtered}\\
v=(4-\partial_x^2)^{-1}u,\qquad
\mathcal A=2v+v_{xx}-3v_x,\qquad
\mathcal B=-(2v+v_{xx}+3v_x).\label{eq:zheng-combinations}
\end{gather}
We deduce that
\begin{equation} 
 \begin{aligned}
 \frac{\mathrm d}{\mathrm dt}[e^q\mathcal A(t,q)]
 &=e^q\left[\frac12\mathcal A^2
                  +\frac32K_-(\mathcal A^2)\right]\ge0,\\
 \frac{\mathrm d}{\mathrm dt}[e^{-q}\mathcal B(t,q)]
 &=e^{-q}\left[\frac12\mathcal B^2
                  +\frac32K_+(\mathcal B^2)\right]\ge0,
 \end{aligned}
\end{equation}
which implies that 
$e^q\mathcal A(t,q)$ and $e^{-q}\mathcal B(t,q)$ are nondecreasing
along each characteristic.

Secondly, we establish the sign relations
\begin{equation}\label{relation} 
 \mathcal A(t,x)>0\Longrightarrow u(t,x)>0,
 \qquad
 \mathcal B(t,x)>0\Longrightarrow u(t,x)<0.
\end{equation}

Thirdly, let $I$ be an open interval where $m_0>0$.
Then $m$ stays positive on $q(t,I)$. We establish a key lemma to show
that $\mathcal A(t,q(t,a))$ becomes positive for some $t\ge0$ and
$a\in I$. An analogous statement holds for $\mathcal B$ on an interval
where $m_0<0$.
 
Finally, the weighted monotonicity gives a common time $T$
and labels $\alpha<\beta$ at which both quantities are positive.
Combining this with \eqref{relation}, we obtain
\[
u(T,q(T,\alpha))>0>u(T,q(T,\beta)).
\]
Since $\alpha<\beta$ and the solution is assumed global, the bound
$u_x\ge-|u|$ gives
\[
u(T,q(T,\beta))\ge e^{q(T,\alpha)-q(T,\beta)}
u(T,q(T,\alpha))>0,
\]
a contradiction.

The converse follows from the global-existence results of Liu and Yin \cite[Theorems~3.1 and~3.2]{LiuYin2006}. This yields Theorem~\ref{thm:momentum}.

\medskip
\noindent The rest of the paper is organized as follows.

\noindent\textbf{Section~\ref{sec:DP-preliminaries}.}
We collect the local theory and the characteristic and momentum identities.

\noindent\textbf{Section~\ref{sec:DP-pointwise}.}
We establish the sharp pointwise criterion and derive the finite-time
characterization.

\noindent\textbf{Section~\ref{mom:section}.}
We prove the momentum characterization of blow-up.

\section{Preliminaries}\label{sec:DP-preliminaries}

We collect the local theory and the identities used in the proofs.
The characteristic and pressure estimates hold for $H^s$ solutions with
$s>3/2$, whereas the pointwise momentum identities are stated for $s>5/2$.

\subsection{Local theory, conserved energy, and pressure}

\begin{lemma}[{\cite[Theorems~2.2, 2.3 and~3.1; Remark~2.1]{Yin2003}}]\label{lem:local-theory}
Let $u_0\in H^s(\R)$ with $s>3/2$. The Cauchy problem \eqref{eq:dp-local}
has a unique maximal strong solution
\[
u\in C([0,T^*);H^s(\R))
\cap C^1([0,T^*);H^{s-1}(\R)),\qquad T^*\in(0,\infty].
\]
The solution depends continuously on the initial data on every compact
subinterval of $[0,T^*)$. If $u_0\in H^{s_1}(\R)$ for some $s_1>s$,
then this higher regularity persists on the same maximal interval.
If $T^*<\infty$, then
\begin{equation}\label{eq:continuation}
\liminf_{t\uparrow T^*}\inf_{x\in\R}u_x(t,x)=-\infty.
\end{equation}
\end{lemma}

\begin{lemma}[{\cite[Lemma~3.1]{LiuYin2006}}]\label{lem:energy}
Let $u_0\in H^s(\R)$ with $s>3/2$, and let $u$ be the corresponding
maximal strong solution. The quadratic energy is conserved:
\begin{equation}\label{eq:energy}
E_0=E(t):=\int_{\R}
u(t,x)(1-\partial_x^2)(4-\partial_x^2)^{-1}u(t,x)\dd x.
\end{equation}
Moreover,
\begin{equation}\label{eq:energy-L2}
\frac14\|u(t)\|_2^2\le E_0\le\|u(t)\|_2^2.
\end{equation}
\end{lemma}

For $\rho>0$ and $f\in C_b(\R)$, recall the one-sided convolution
operators from \eqref{eq:intro-K}
\begin{equation}\label{eq:K}
\begin{aligned}
K_{\rho,-}f(x)&=\int_{-\infty}^x e^{\rho(y-x)}f(y)\dd y,\\
K_{\rho,+}f(x)&=\int_x^\infty e^{\rho(x-y)}f(y)\dd y.
\end{aligned}
\end{equation}
Here $C_b(\R)$ denotes the bounded continuous functions, and we write
$K_-=K_{1,-}$ and $K_+=K_{1,+}$.

\begin{lemma}[One-sided convolution and pressure estimates]\label{lem:K-pressure}
For $\rho>0$ and $f\in C_b(\R)$, the functions $K_{\rho,\pm}f$
belong to $C_b^1(\R)$ and satisfy
\begin{equation}\label{eq:K-identities}
\partial_xK_{\rho,-}f=f-\rho K_{\rho,-}f,\qquad
\partial_xK_{\rho,+}f=\rho K_{\rho,+}f-f,
\end{equation}
and
\begin{equation}\label{eq:K-norm}
\|K_{\rho,\pm}f\|_\infty\le\rho^{-1}\|f\|_\infty.
\end{equation}
These operators preserve nonnegativity. If $u_0\in H^s(\R)$ with
$s>3/2$ and $u$ is the corresponding maximal strong solution, then
the pressure $\Pi=\frac32p*u^2$ satisfies
\begin{equation}\label{eq:pressure}
\Pi-\Pi_x=\frac32K_-(u^2),\qquad
\Pi+\Pi_x=\frac32K_+(u^2),\qquad
\Pi_{xx}=\Pi-\frac32u^2.
\end{equation}
Consequently,
\begin{equation}\label{eq:pressure-bound}
0\le |\Pi_x|\le\Pi\le\frac34\|u(t)\|_2^2\le3E_0,
\end{equation}
and
\begin{equation}\label{eq:pressure-Linfty}
\Pi\le\frac32\|u(t)\|_\infty^2,
\qquad
0\le\Pi\pm\Pi_x\le\frac32\|u(t)\|_\infty^2.
\end{equation}
\end{lemma}

\begin{proof}
For the pressure representation \eqref{eq:dp} and its bound in terms
of $\|u\|_2$, see
\cite[(2.4)]{Yin2003} and \cite[proof of Lemma~3.2]{LiuYin2006}.
In the present notation,
\[
\begin{aligned}
\partial_x\left(e^{-\rho x}\int_{-\infty}^x e^{\rho y}f(y)\dd y\right)
&=f(x)-\rho K_{\rho,-}f(x),\\[0.3em]
\partial_x\left(e^{\rho x}\int_x^\infty e^{-\rho y}f(y)\dd y\right)
&=\rho K_{\rho,+}f(x)-f(x),
\end{aligned}
\]
\[
|K_{\rho,\pm}f(x)|
\le\|f\|_\infty\int_0^\infty e^{-\rho r}\dd r
=\rho^{-1}\|f\|_\infty.
\]
The kernels are nonnegative, and the Green kernel satisfies
\begin{equation}\label{eq:green-kernel}
p-p_{xx}=\delta_0.
\end{equation}
\[
\begin{aligned}
\Pi&=\frac34\bigl(K_-(u^2)+K_+(u^2)\bigr),\\[0.3em]
\Pi_x&=\frac34\bigl(K_+(u^2)-K_-(u^2)\bigr).
\end{aligned}
\]
Using \eqref{eq:K-norm} and \eqref{eq:energy-L2}, we obtain
\[
\begin{gathered}
\Pi\mp\Pi_x=\frac32K_{\mp}(u^2),\qquad
\Pi_{xx}=\Pi-\frac32u^2,\\[0.3em]
|\Pi_x|\le\Pi
\le\frac32\|p\|_\infty\|u(t)\|_2^2
=\frac34\|u(t)\|_2^2\le3E_0,\\[0.3em]
\Pi\le\frac32\|p\|_1\|u(t)\|_\infty^2
=\frac32\|u(t)\|_\infty^2,\\[0.3em]
0\le\Pi\mp\Pi_x
=\frac32K_{\mp}(u^2)
\le\frac32\|u(t)\|_\infty^2.
\end{gathered}
\]
\end{proof}

\subsection{Characteristic identities}

Recall the characteristic flow and the Eulerian quantities $A,B$
from \eqref{eq:intro-flow} and \eqref{eq:intro-AB}
\begin{align}
q_t(t,a)&=u(t,q(t,a)),\qquad q(0,a)=a,\label{eq:flow}\\
A(t,x)&=u(t,x)-u_x(t,x),\qquad
B(t,x)=-u(t,x)-u_x(t,x).\label{eq:AB}
\end{align}
Here $a$ denotes a Lagrangian label and $x$ the Eulerian variable.
All convolutions act on the Eulerian variable. We have
\begin{equation}\label{eq:AB-algebra}
A-B=2u,\quad A+B=-2u_x,\quad
AB=u_x^2-u^2,\quad \min\{A,B\}=-u_x-|u|.
\end{equation}
\begin{lemma}[Characteristic identities]\label{lem:characteristic}
Let $u_0\in H^s(\R)$ with $s>3/2$, and let $u$ be the corresponding
maximal strong solution. The flow $q$ satisfies
\begin{equation}\label{eq:qa}
q_a(t,a)=\exp\left(\int_0^t u_x(\tau,q(\tau,a))\dd\tau\right)>0.
\end{equation}
On every compact time interval in $[0,T^*)$, $|q(t,a)-a|$ is bounded
uniformly in $a$. Hence $q(t,\cdot)$ is an increasing $C^1$
diffeomorphism of $\R$. Put
\begin{equation}\label{eq:UW}
U(t,a)=u(t,q(t,a)),\qquad W(t,a)=u_x(t,q(t,a)).
\end{equation}

For each fixed label $a$, the functions $U(t,a)$ and $W(t,a)$ are
continuously differentiable in time, and
\begin{equation}\label{sp:characteristic-identities}
\begin{aligned}
\partial_tU(t,a)&=-\Pi_x(t,q(t,a)),\\
\partial_tW(t,a)&=-W(t,a)^2+\frac32U(t,a)^2-\Pi(t,q(t,a)).
\end{aligned}
\end{equation}
Moreover,
\begin{equation}\label{sp:AB-flow}
\begin{aligned}
\dt A(t,q(t,a))
&=W^2-\frac32U^2+\frac32K_-(u^2)(t,q(t,a)),\\
\dt B(t,q(t,a))
&=W^2-\frac32U^2+\frac32K_+(u^2)(t,q(t,a)).
\end{aligned}
\end{equation}

\end{lemma}

\begin{proof}
For the classical characteristic formulas, see
\cite[Lemma~2.3 and proof of Lemma~3.2]{LiuYin2006}.
Fix $T<T^*$. Lemma~\ref{lem:local-theory} and Sobolev embedding give
\begin{equation}\label{eq:strong-pointwise-regularity}
u\in C([0,T];C_b^1(\R)),\qquad u_t\in C([0,T];C_b(\R)).
\end{equation}
By \eqref{eq:pressure}, $\Pi_{xx}=\Pi-\frac32u^2$ also belongs to
$C([0,T];C_b(\R))$.
By \eqref{eq:flow},
\[
\partial_tq_a=Wq_a,\qquad q_a(0,a)=1,
\qquad
q_a(t,a)=\exp\left(\int_0^t W(\tau,a)\dd\tau\right)>0,
\]
\[
|q(t,a)-a|\le\int_0^t\|u(\tau)\|_\infty\dd\tau,
\qquad \lim_{a\to\pm\infty}q(t,a)=\pm\infty.
\]
The chain rule along a fixed characteristic and \eqref{eq:dp} give
\begin{equation}\label{sp:characteristic-integral}
U(t,a)=u_0(a)-\int_0^t\Pi_x(\tau,q(\tau,a))\dd\tau.
\end{equation}
Differentiating in $a$ gives
\begin{equation}\label{sp:label-differentiation}
W(t,a)q_a(t,a)
=u_0'(a)-\int_0^t
\Pi_{xx}(\tau,q(\tau,a))q_a(\tau,a)\dd\tau.
\end{equation}
By \eqref{eq:qa} and \eqref{sp:label-differentiation},
\[
\begin{gathered}
W(\cdot,a)q_a(\cdot,a)\in C^1([0,T]),\quad
q_a(\cdot,a)\in C^1([0,T]),\quad q_a>0\\[0.3em]
\Longrightarrow\quad W(\cdot,a)\in C^1([0,T]).
\end{gathered}
\]
Using \eqref{eq:pressure} in the time derivatives of
\eqref{sp:characteristic-integral} and \eqref{sp:label-differentiation}
gives
\[
\begin{aligned}
\partial_tU&=-\Pi_x(t,q),\\[0.3em]
\partial_tW
&=\frac{\partial_t(Wq_a)-W\partial_tq_a}{q_a}
=-\Pi_{xx}(t,q)-W^2\\[0.3em]
&=-W^2+\frac32U^2-\Pi(t,q).
\end{aligned}
\]
By \eqref{eq:AB} and \eqref{eq:UW},
\[
\begin{aligned}
\frac{\mathrm d}{\mathrm dt}A(t,q)
&=\partial_tU-\partial_tW
=W^2-\frac32U^2+(\Pi-\Pi_x)(t,q),\\[0.3em]
\frac{\mathrm d}{\mathrm dt}B(t,q)
&=-\partial_tU-\partial_tW
=W^2-\frac32U^2+(\Pi+\Pi_x)(t,q).
\end{aligned}
\]
\end{proof}

Integrating the first identity in \eqref{sp:characteristic-identities}
and using \eqref{eq:pressure-bound} and the surjectivity in
Lemma~\ref{lem:characteristic}, we obtain
\begin{equation}\label{eq:finite-amplitude}
\|u(t)\|_\infty\le\|u_0\|_\infty+3E_0t.
\end{equation}

With $q=q(t,a)$, the identities \eqref{sp:AB-flow} can also be written as
\begin{equation}\label{eq:raw-AB-product}
\begin{aligned}
\frac{\mathrm d}{\mathrm dt}A(t,q)
&=\left[AB+\frac12K_-(A^2)-\frac12K_-(AB)\right](t,q),\\
\frac{\mathrm d}{\mathrm dt}B(t,q)
&=\left[AB+\frac12K_+(B^2)-\frac12K_+(AB)\right](t,q).
\end{aligned}
\end{equation}
Indeed, \eqref{eq:AB} and integration by parts in \eqref{eq:K} give
\begin{equation}\label{eq:raw-AB-square-integrals}
\begin{aligned}
\frac12K_-(A^2)
&=K_-(u^2)+\frac12K_-(u_x^2)-\frac12u^2,\\
\frac12K_+(B^2)
&=K_+(u^2)+\frac12K_+(u_x^2)-\frac12u^2.
\end{aligned}
\end{equation}
Combining these identities with \eqref{eq:AB-algebra} and
\eqref{sp:AB-flow} yields \eqref{eq:raw-AB-product}.

\begin{remark}[Comparison with the Camassa--Holm equation]
For the Camassa--Holm equation, the same quantities satisfy
\begin{equation}\label{eq:CH-raw-AB-comparison}
\begin{aligned}
\frac{\mathrm d}{\mathrm dt}A(t,q)
&=\frac12AB+\frac12K_-(A^2),\\
\frac{\mathrm d}{\mathrm dt}B(t,q)
&=\frac12AB+\frac12K_+(B^2);
\end{aligned}
\end{equation}
see \cite[Section~2, formulas (3) and (4)]{McKean1998}.
In particular, the corresponding weighted quantities satisfy
\[
\begin{aligned}
\frac{\mathrm d}{\mathrm dt}\bigl[e^qA(t,q)\bigr]
&=\frac12e^q\bigl[A^2+K_-(A^2)\bigr]\ge0,\\
\frac{\mathrm d}{\mathrm dt}\bigl[e^{-q}B(t,q)\bigr]
&=\frac12e^{-q}\bigl[B^2+K_+(B^2)\bigr]\ge0.
\end{aligned}
\]
For the DP equation, \eqref{eq:raw-AB-product} instead gives
\begin{equation}\label{eq:DP-raw-AB-weighted}
\begin{aligned}
\frac{\mathrm d}{\mathrm dt}\bigl[e^{2q}A(t,q)\bigr]
&=e^{2q}\left[A^2+\frac12K_-(A^2)-\frac12K_-(AB)\right],\\
\frac{\mathrm d}{\mathrm dt}\bigl[e^{-2q}B(t,q)\bigr]
&=e^{-2q}\left[B^2+\frac12K_+(B^2)-\frac12K_+(AB)\right].
\end{aligned}
\end{equation}
The additional terms involving $K_\pm(AB)$ have no prescribed sign.
These identities do not ensure that $A$ and $B$ remain positive along
a fixed characteristic.
\end{remark}

\subsection{The spatial maximum}

Recall from \eqref{eq:intro-J} that
\begin{equation}\label{eq:J}
J(t):=\sup_{x\in\R}\bigl(\min\{A(t,x),B(t,x)\}\bigr)_+
=\sup_{x\in\R}\bigl(-u_x(t,x)-|u(t,x)|\bigr)_+,
\end{equation}
where $r_+=\max\{r,0\}$. By \eqref{eq:J-slope} and
\eqref{eq:finite-amplitude}, at a finite endpoint, $J(t)\to+\infty$
if and only if $\inf_{x\in\R}u_x(t,x)\to-\infty$. By Sobolev embedding,
$u(t,\cdot),u_x(t,\cdot)\in C_0(\R)$, so the supremum is attained
whenever $J(t)>0$.

Set
\begin{equation}\label{sp:F}
F(t,a):=-W(t,a)-|U(t,a)|
=\min\{A(t,q(t,a)),B(t,q(t,a))\}.
\end{equation}
Since $q(t,\cdot)$ is onto by Lemma~\ref{lem:characteristic},
\eqref{eq:J} gives
\begin{equation}\label{eq:J-label}
J(t)=\sup_{a\in\R}F(t,a)_+,
\qquad u_x(t,x)\ge-|u(t,x)|-J(t)\quad(x\in\R).
\end{equation}
\begin{lemma}[Time regularity of the spatial maximum]\label{lem:J-regularity}
Let $u_0\in H^s(\R)$ with $s>3/2$, and let $u$ be the corresponding
maximal strong solution. Then $J$ is locally Lipschitz on $[0,T^*)$.
In particular, it is absolutely continuous on every $[0,T]$ with
$T<T^*$ and differentiable almost everywhere in $(0,T^*)$.
\end{lemma}

\begin{proof}
Fix $T<T^*$.
By \eqref{eq:strong-pointwise-regularity}, $U$ and $W$ are uniformly
bounded on $[0,T]\times\R$. Equations
\eqref{sp:characteristic-identities} and \eqref{eq:pressure-Linfty}
then bound their time derivatives on this set.
Hence there is a finite constant $C_T$ such that
\begin{equation}\label{sp:label-lipschitz}
|U(t,a)-U(r,a)|+|W(t,a)-W(r,a)|
\le C_T|t-r|,
\qquad t,r\in[0,T],\quad a\in\R.
\end{equation}
Since the absolute value and positive-part maps are $1$-Lipschitz,
\eqref{eq:J-label} gives
\begin{equation}\label{sp:J-lipschitz}
|J(t)-J(r)|
\le\sup_{a\in\R}|F(t,a)_+-F(r,a)_+|
\le C_T|t-r|.
\end{equation}
\end{proof}

\subsection{Momentum identities}

\begin{lemma}[Momentum transport and its consequences]\label{lem:momentum-identities}
Let $u_0\in H^s(\R)$ with $s>5/2$, and let $u$ be the corresponding
maximal strong solution. Then $m=u-u_{xx}$ belongs to
$C([0,T^*);H^{s-2}(\R))$, and $m(t,\cdot)$ is continuous and bounded.
The momentum transport identity, stated for $s\ge3$ in
\cite[Lemma~2.4]{LiuYin2006}, holds in the form
\begin{equation}\label{eq:momentum-transport}
m(t,q(t,a))q_a(t,a)^3=m_0(a).
\end{equation}
The first-order spatial identities
\begin{equation}\label{eq:momentum-spatial}
A_x+A=m,\qquad B_x-B=m
\end{equation}
also yield the one-sided momentum representations
\cite[formula (3.3)]{LiuYin2007}
\begin{equation}\label{eq:momentum-convolution}
A=K_-m,\qquad B=-K_+m.
\end{equation}
For every finite label interval $[a,b]$,
\begin{equation}\label{eq:local-momentum-invariant}
\int_{q(t,a)}^{q(t,b)}|m(t,x)|^{1/3}\dd x
=\int_a^b|m_0(r)|^{1/3}\dd r.
\end{equation}
For labels $a<b$, write $L_{ab}(t)=q(t,b)-q(t,a)$. Then
\begin{equation}\label{eq:interval-log-derivative}
\frac{\mathrm d}{\mathrm dt}\log L_{ab}(t)
=\frac{1}{L_{ab}(t)}\int_{q(t,a)}^{q(t,b)}u_x(t,x)\dd x.
\end{equation}
\end{lemma}

\begin{proof}
Fix $T<T^*$. For smooth data, \eqref{eq:momentum} and \eqref{eq:qa}
give
\[
\frac{\mathrm d}{\mathrm dt}\bigl[m(t,q)q_a^3\bigr]
=\bigl(m_t+um_x+3u_xm\bigr)(t,q)q_a^3=0.
\]
Choose $u_{0,n}\in C_c^\infty(\R)$ with $u_{0,n}\to u_0$ in $H^s$,
and denote the corresponding quantities by $u_n,q_n,m_n$.
By Lemma~\ref{lem:local-theory}, the solutions $u_n$ remain smooth on
$[0,T]$ for all sufficiently large $n$ and converge to $u$ in
$C([0,T];H^s(\R))$. Sobolev embedding, Gr\"onwall's inequality applied
to \eqref{eq:flow}, and the formula \eqref{eq:qa} for the label
derivatives then give
\[
\begin{alignedat}{2}
u_n&\longrightarrow u
&\qquad&\text{in }C([0,T];H^s(\R))\hphantom{\text{in }}\hookrightarrow C([0,T];C_b^2(\R)),\\[0.3em]
m_n&\longrightarrow m
&\qquad&\text{in }C([0,T];C_b(\R)),\\[0.3em]
(q_n,\partial_aq_n)&\longrightarrow(q,q_a)
&\qquad&\text{locally uniformly on }[0,T]\times\R.
\end{alignedat}
\]
Therefore,
\[
m_n(t,q_n)(\partial_aq_n)^3=m_{0,n}
\quad\xrightarrow{\ n\to\infty\ }\quad
m(t,q)q_a^3=m_0.
\]
By \eqref{eq:AB} and \eqref{eq:momentum},
\[
\begin{aligned}
\partial_x(e^xA)&=e^x(A_x+A)=e^x(u-u_{xx})=e^xm,\\[0.3em]
\partial_x(e^{-x}B)&=e^{-x}(B_x-B)=e^{-x}(u-u_{xx})=e^{-x}m,
\end{aligned}
\]
\[
\lim_{x\to-\infty}e^xA(t,x)=0,\qquad
\lim_{x\to\infty}e^{-x}B(t,x)=0,
\]
\[
A(t,x)=e^{-x}\int_{-\infty}^x e^ym(t,y)\dd y=K_-m(t,x),
\]
\[
B(t,x)=-e^x\int_x^\infty e^{-y}m(t,y)\dd y=-K_+m(t,x).
\]
Using \eqref{eq:momentum-transport} in the change of variables gives
\[
\begin{aligned}
\int_{q(t,a)}^{q(t,b)}|m(t,x)|^{1/3}\dd x
&=\int_a^b|m(t,q(t,r))|^{1/3}q_r(t,r)\dd r\\
&=\int_a^b|m_0(r)|^{1/3}\dd r.
\end{aligned}
\]
The flow equation \eqref{eq:flow} gives
\[
 L_{ab}'(t)=u(t,q(t,b))-u(t,q(t,a))
 =\int_{q(t,a)}^{q(t,b)}u_x(t,x)\dd x.
\]
Dividing by $L_{ab}(t)>0$ yields \eqref{eq:interval-log-derivative}.
\end{proof}

\subsection{One-sided auxiliary quantities}

Let $u_0\in H^s(\R)$ with $s>3/2$, and let $u$ be the corresponding
maximal strong solution. Recall the quantities from
\eqref{eq:intro-filtered}
\begin{equation}\label{mom:filtered-definitions}
 \begin{aligned}
 \mathcal A(t,x)
 &=K_{2,-}A(t,x)
 =\int_{-\infty}^xe^{2(y-x)}(u-u_x)(t,y)\dd y,\\
 \mathcal B(t,x)
 &=K_{2,+}B(t,x)
 =\int_x^\infty e^{2(x-y)}(-u-u_x)(t,y)\dd y.
 \end{aligned}
\end{equation}
The bound \eqref{eq:K-norm} gives
\begin{equation}\label{eq:filtered-norms}
\|\mathcal A(t)\|_\infty\le\frac12\|A(t)\|_\infty,
\qquad
\|\mathcal B(t)\|_\infty\le\frac12\|B(t)\|_\infty.
\end{equation}
By \eqref{eq:AB} and \eqref{eq:K-identities},
\begin{equation}\label{mom:filtered-algebra}
 \begin{gathered}
 \mathcal A=3K_{2,-}u-u,\qquad
 \mathcal B=u-3K_{2,+}u,\\
 (\partial_x+2)\mathcal A=A,\qquad
 (\partial_x-2)\mathcal B=-B,\\
 m=(\partial_x+1)(\partial_x+2)\mathcal A
   =-(\partial_x-1)(\partial_x-2)\mathcal B.
 \end{gathered}
\end{equation}
The identities involving $m$ hold distributionally for $s>3/2$ and
pointwise for $s>5/2$.
For $v$ in \eqref{eq:zheng-combinations}, the kernel of
$(4-\partial_x^2)^{-1}$ is $\frac14e^{-2|x|}$. Thus
$K_{2,-}u=2v-v_x$ and $K_{2,+}u=2v+v_x$ by \eqref{eq:K-identities}.
Using $u=4v-v_{xx}$ in \eqref{mom:filtered-algebra} gives
\eqref{eq:zheng-combinations}; see \cite[Section~3, formulas (3.8) and
(3.9)]{Zheng2025} for these combinations.

The evenness of $p$ in \eqref{eq:dp} gives
$\widetilde\Pi(t,x)=\Pi(t,-x)$ for $\widetilde u(t,x)=-u(t,-x)$,
and substitution into \eqref{eq:dp} shows that $\widetilde u$ solves
the same equation. By \eqref{eq:momentum}, uniqueness for
\eqref{eq:flow}, and \eqref{mom:filtered-definitions},
\begin{equation}\label{eq:reflection}
\begin{gathered}
\widetilde\Pi(t,x)=\Pi(t,-x),\qquad
\widetilde m(t,x)=-m(t,-x),\qquad
\widetilde q(t,a)=-q(t,-a),\\
\widetilde A(t,x)=B(t,-x),\qquad
\widetilde B(t,x)=A(t,-x),\\
\widetilde{\mathcal A}(t,x)=\mathcal B(t,-x),\qquad
\widetilde{\mathcal B}(t,x)=\mathcal A(t,-x).
\end{gathered}
\end{equation}

\begin{lemma}[Identities for the auxiliary quantities]\label{mom:filtered-identities}
Let $u_0\in H^s(\R)$ with $s>3/2$, and let $u$ be the corresponding
{ strong solution}. Along every characteristic,
\begin{equation}\label{mom:filtered-evolution}
 \begin{aligned}
 \frac{\mathrm d}{\mathrm dt}\mathcal A(t,q)
 &=\left(\frac12\mathcal A-u\right)\mathcal A
   +\frac32K_-(\mathcal A^2),\\
 \frac{\mathrm d}{\mathrm dt}\mathcal B(t,q)
 &=\left(\frac12\mathcal B+u\right)\mathcal B
   +\frac32K_+(\mathcal B^2),
 \end{aligned}
\end{equation}
where all terms on the right-hand sides are evaluated at $(t,q(t,a))$.
Consequently,
\begin{equation}\label{mom:weighted-monotonicity}
 \begin{aligned}
 \frac{\mathrm d}{\mathrm dt}[e^q\mathcal A(t,q)]
 &=e^q\left[\frac12\mathcal A^2
                  +\frac32K_-(\mathcal A^2)\right]\ge0,\\
 \frac{\mathrm d}{\mathrm dt}[e^{-q}\mathcal B(t,q)]
 &=e^{-q}\left[\frac12\mathcal B^2
                  +\frac32K_+(\mathcal B^2)\right]\ge0.
 \end{aligned}
\end{equation}
\end{lemma}

\begin{proof}
Put $h=K_{2,-}u$. By \eqref{eq:strong-pointwise-regularity} and
\eqref{eq:K-norm}, $h_t=K_{2,-}u_t$. Equations
\eqref{eq:K-identities} and \eqref{mom:filtered-algebra} give
\[
 h_x=u-2h,\qquad \mathcal A=3h-u=h-h_x,
\]
so $\mathcal A$ is continuously differentiable in both variables.
By \eqref{eq:dp} and integration by parts in \eqref{eq:K},
\[
 h_t=-\frac12u^2+K_{2,-}(u^2)-\Pi+2K_{2,-}\Pi.
\]
Using \eqref{eq:pressure}, we obtain
\[
 3K_{2,-}(u^2)-3\Pi+6K_{2,-}\Pi+\Pi_x
 =\frac32K_-(u^2).
\]
Both sides are bounded, and \eqref{eq:K-identities} and
\eqref{eq:pressure} give the same expression after applying
$\partial_x+2$. The only bounded solution of the resulting
homogeneous equation for the difference is zero. Hence
\[
 (\partial_t+u\partial_x)\mathcal A
 =-2u\mathcal A-\frac12u^2+\frac32K_-(u^2).
\]
By \eqref{eq:K} and integration by parts,
\[
 \begin{aligned}
 u^2-\mathcal A^2
 &=(h_x+2h)^2-(h-h_x)^2
 =3(h^2)_x+3h^2,\\
 K_-(u^2)&=K_-(\mathcal A^2)+3h^2.
 \end{aligned}
\]
Substitution of $h=(u+\mathcal A)/3$ proves the first identity in
\eqref{mom:filtered-evolution}.
The reflection identities \eqref{eq:reflection} give the second
identity. Multiplication by $e^q$ and $e^{-q}$ and use of
\eqref{eq:flow} give \eqref{mom:weighted-monotonicity}.

\end{proof}

\section{The sharp pointwise criterion and the finite-time characterization}\label{sec:DP-pointwise}

We first establish the differential inequality for $J$ and prove the
pointwise criterion. We then verify the optimality of its coefficient
and derive the finite-time characterization.

\subsection{The pointwise criterion}

\begin{proof}[Proof of Theorem~\ref{thm:single}]
The proof is divided into four steps.
\proofstep{Step 1. A maximizing label at a fixed time.}
By Lemma~\ref{lem:J-regularity}, $J$ is locally Lipschitz and $J(0)>0$ so that $J(t)\not\equiv0$ in $(0,T^*)$.
Fix a time $t\in(0,T^*)$ at which $J$ is differentiable.
If $J(t)=0$, then $J'(t)=0$ because $J\ge0$.
Suppose that $J(t)>0$. The supremum in \eqref{eq:J} is attained,
and $q(t,\cdot)$ is a diffeomorphism by Lemma~\ref{lem:characteristic}.
Choose $a$ as the inverse image of a spatial maximizer under
$q(t,\cdot)$. Then $F(t,a)=J(t)$ by \eqref{eq:J-label}.
With the notation of
\eqref{eq:UW} and \eqref{sp:F}, write
\begin{equation}\label{sp:active-slope}
J=J(t),\qquad U=U(t,a),\qquad
W=W(t,a)=-|U|-J.
\end{equation}
Keep this label fixed as time varies. By \eqref{eq:J-label},
\[
J(t+h)\ge F(t+h,a)_+\ge F(t+h,a)
\qquad(0<h<T^*-t).
\]
Since $J(t)=F(t,a)$, comparison of the right time derivatives gives
\begin{equation}\label{sp:envelope-comparison}
J'(t)\ge\partial_t^+F(t,a),
\end{equation}
where $\partial_t^+$ denotes the right time derivative at fixed $a$.
This derivative exists by \eqref{sp:F} and the time regularity of
$U,W$ in Lemma~\ref{lem:characteristic}.

\proofstep{Step 2. The case \(U\ge0\).}
Suppose first that $U>0$. We first derive a pointwise lower bound for
$K_+(u^2)(t,q(t,a))$. For the spatial estimate, keep $t$ fixed and set
\[
f(y)=u(t,q(t,a)+y)\quad(y\ge0),\qquad
L=\log\frac{U+J}{J}>0.
\]
On the interval where $f$ stays positive starting from $f(0)=U$,
\eqref{eq:J-label} gives
\[
f'(y)+f(y)\ge-J,\qquad
\frac{\mathrm d}{\mathrm dy}\bigl[e^y(f(y)+J)\bigr]\ge0.
\]
Hence $f(y)\ge(U+J)e^{-y}-J$. This lower bound is strictly positive for $0\le y<L$ by the definition of $L$.
By continuity,
\begin{equation}\label{sp:barrier-interval}
f(y)\ge(U+J)e^{-y}-J\ge0\qquad(0\le y\le L).
\end{equation}
Using \eqref{sp:barrier-interval}, the definition \eqref{eq:K} of $K_+$,
and $z=e^{-y}$, we obtain
\begin{align}
K_+(u^2)(t,q(t,a))
&=\int_0^\infty e^{-y}f(y)^2\dd y\nonumber\\
&\ge\int_0^L e^{-y}\bigl((U+J)e^{-y}-J\bigr)^2\dd y\nonumber\\
&=\int_{J/(U+J)}^1\bigl((U+J)z-J\bigr)^2\dd z
=\frac{U^3}{3(U+J)}.
\label{sp:right-pressure}
\end{align}

Now vary time along the fixed label $a$. By continuity of $U(\cdot,a)$
and \eqref{eq:AB}, \eqref{sp:F} yields $F(\tau,a)=B(\tau,q(\tau,a))$
for $\tau$ sufficiently close to $t$. Using \eqref{sp:envelope-comparison},
\eqref{sp:AB-flow}, \eqref{sp:active-slope}, and
\eqref{sp:right-pressure}, we obtain
\begin{align}
J'(t)
&\ge W^2-\frac32U^2+\frac32K_+(u^2)(t,q(t,a))\nonumber\\
&\ge(U+J)^2-\frac32U^2+\frac{U^3}{2(U+J)}\nonumber\\
&=J^2+\frac{UJ(3U+4J)}{2(U+J)}\ge J^2.
\label{sp:positive-growth}
\end{align}

If $U(t,a)=0$, differentiability of $U(\cdot,a)$ gives
\[
\begin{aligned}
\partial_t^+|U(t,a)|
&=\lim_{h\downarrow0}\frac{|U(t+h,a)|-|U(t,a)|}{h}\\
&=\lim_{h\downarrow0}
  \left|\frac{U(t+h,a)-U(t,a)}{h}\right|\\
&=|\partial_tU(t,a)|.
\end{aligned}
\]
By \eqref{sp:F}, \eqref{sp:envelope-comparison},
\eqref{sp:active-slope}, \eqref{sp:characteristic-identities}, and
\eqref{eq:pressure-bound},
\begin{equation}\label{sp:zero-growth}
\begin{aligned}
J'(t)&\ge\partial_t^+F(t,a)\\
&=-\partial_tW(t,a)-\partial_t^+|U(t,a)|\\
&=W(t,a)^2-\frac32U(t,a)^2
  +\Pi(t,q(t,a))-|\Pi_x(t,q(t,a))|\\
&=J^2+\Pi(t,q(t,a))-|\Pi_x(t,q(t,a))|\ge J^2.
\end{aligned}
\end{equation}

\proofstep{Step 3. The case \(U<0\).}
Consider the reflected solution in \eqref{eq:reflection}.
The spatial definition \eqref{eq:J} gives
$\widetilde J(\tau)=J(\tau)$ at every time $\tau$.
By \eqref{eq:reflection} and \eqref{eq:UW},
\[
\widetilde q(t,-a)=-q(t,a),\qquad
\widetilde U(t,-a)=-U>0,\qquad
\widetilde W(t,-a)=W.
\]
Thus $\widetilde F(t,-a)=F(t,a)=\widetilde J(t)$ by \eqref{sp:F}.
Applying Step 2 at the label $-a$ gives $J'(t)\ge J^2$.

\proofstep{Step 4. Growth of \(J\) and finite lifespan.}
Combining the preceding cases gives
\begin{equation}\label{sp:J-growth}
J'(t)\ge J(t)^2\quad\text{for a.e. }t\in(0,T^*).
\end{equation}
By Lemma~\ref{lem:J-regularity}, $J$ is absolutely continuous on
every $[0,T]$ with $T<T^*$. Integrating \eqref{sp:J-growth} shows that
$J$ is nondecreasing on $[0,T^*)$. If $J(t_0)>0$ for some
$t_0\in[0,T^*)$, then
$J(t)\ge J(t_0)>0$ for all $t\in[t_0,T^*)$, and $1/J$ is locally
absolutely continuous there. Integrating
\[
\left(\frac1{J(t)}\right)'=-\frac{J'(t)}{J(t)^2}\le-1
\quad\text{a.e.}
\]
gives
\begin{equation}\label{sp:reciprocal}
0<\frac1{J(t)}\le\frac1{J(t_0)}-(t-t_0),
\qquad t_0\le t<T^*.
\end{equation}
Therefore
\begin{equation}\label{sp:lifespan}
T^*\le t_0+\frac1{J(t_0)}<\infty.
\end{equation}
In particular, if \eqref{eq:single-initial} holds, then \eqref{eq:J} gives
$J(0)\ge-u_0'(x_0)-|u_0(x_0)|>0$, and
\begin{equation}\label{sp:initial-point-lifespan}
T^*\le\frac1{J(0)}\le\frac1{-u_0'(x_0)-|u_0(x_0)|}.
\end{equation}
\end{proof}

\subsection{Sharpness of the pointwise coefficient}

\begin{proposition}[Sharpness of the pointwise criterion]\label{prop:sharp}
There exist $u_0\in\bigcap_{r\ge0}H^r(\R)$ and $x_0\in\R$ such that
the corresponding strong solution of \eqref{eq:dp} is global and
\[
u_0'(x_0)=-|u_0(x_0)|<0.
\]
In particular, the same datum satisfies
\[
u_0'(x_0)<-\beta|u_0(x_0)|\qquad\text{for every }\beta<1.
\]
Thus the coefficient $1$ in Theorem~\ref{thm:single} cannot be reduced.
\end{proposition}

\begin{proof}[Proof of Proposition~\ref{prop:sharp}]
Choose any nonzero function
\[
m_0\in C_c^\infty(\R),\qquad m_0\ge0,
\]
and set $u_0=p*m_0$.
By the Green kernel identity \eqref{eq:green-kernel},
$u_0-u_0''=m_0$ and $u_0\in\bigcap_{r\ge0}H^r(\R)$.
By \eqref{eq:momentum-transport} and the surjectivity of the flow
in Lemma~\ref{lem:characteristic}, $m(t,\cdot)\ge0$ throughout
the maximal interval of existence.
By \eqref{eq:green-kernel}, $u=p*m$. The kernel inequality
$|p'|\le p$ for $p$ in \eqref{eq:dp} then gives
\begin{equation}\label{sp:sharp-slope}
u(t,x)\ge0,\qquad |u_x(t,x)|\le u(t,x).
\end{equation}
If $T^*<\infty$, \eqref{sp:sharp-slope} and \eqref{eq:finite-amplitude}
would give
\[
u_x(t,x)\ge-u(t,x)
\ge-\|u_0\|_\infty-3E_0T^*
\qquad(0\le t<T^*),
\]
contradicting \eqref{eq:continuation}. Hence the solution is global.

Choose $R$ such that $\operatorname{supp}m_0\subset[-R,R]$.
With $p$ as in \eqref{eq:dp}, for every $x>R$,
\begin{equation}\label{sp:sharp-tail}
u_0(x)=\frac12e^{-x}\int_{\R}e^ym_0(y)\dd y>0,
\qquad u_0'(x)=-u_0(x).
\end{equation}
Choose any $x_0>R$. By \eqref{sp:sharp-tail},
$u_0'(x_0)=-|u_0(x_0)|<0$, so the same datum satisfies the strict
condition with every $\beta<1$.
\end{proof}

\subsection{The finite-time characterization}

\begin{proof}[Proof of Theorem~\ref{thm:entry}]
If the condition on the right of \eqref{eq:entry-iff} holds,
then $J(t_1)>0$ by \eqref{eq:J}. Taking $t_0=t_1$ in
\eqref{sp:lifespan} gives $T^*<\infty$.

Conversely, suppose that $T^*<\infty$. By \eqref{eq:finite-amplitude},
\[
\sup_{0\le t<T^*}\|u(t)\|_\infty
\le M_*:=\|u_0\|_\infty+3E_0T^*<\infty.
\]
The continuation criterion \eqref{eq:continuation} and
\eqref{eq:J-slope} yield
\[
-\inf_{x\in\R}u_x(t,x)-M_*\le J(t),\qquad
\limsup_{t\uparrow T^*}J(t)=+\infty.
\]
In particular, $J(t_1)>0$ for some $t_1\in[0,T^*)$, and \eqref{eq:J}
gives $x_1\in\R$ with $u_x(t_1,x_1)<-|u(t_1,x_1)|$.
Moreover, since $J$ is nondecreasing by \eqref{sp:J-growth}
and Lemma~\ref{lem:J-regularity}, $J(t)\to+\infty$ as $t\uparrow T^*$.
It follows again from \eqref{eq:J-slope} that
\[
\inf_{x\in\R}u_x(t,x)\le-J(t)\longrightarrow-\infty
\qquad(t\uparrow T^*).
\]
This and the amplitude bound give \eqref{eq:breaking}.
\end{proof}

\section{The momentum characterization of blow-up}
\label{mom:section}

Throughout this section, $s>5/2$ and $u_0\in H^s(\R)$.
We first establish uniform bounds for global solutions and derive a
local sign property on an interval of positive momentum. The weighted
identities then allow us to use the signs on two ordered intervals at
a common time to prove Theorem~\ref{thm:momentum}.

\subsection{Uniform bounds for global solutions}

\begin{lemma}[Uniform bounds for global solutions]\label{mom:global-bounds}
If $T^*=\infty$, then
\begin{equation}\label{mom:global-estimates}
 \begin{gathered}
 u_x(t,x)\ge-|u(t,x)|,\qquad |u(t,x)|\le M:=\sqrt{8E_0},\\
 -M\le u_x(t,x)\le C,
 \qquad C:=\max\{\sqrt{3/2}\,M,\|u_0'\|_{L^\infty}\}.
 \end{gathered}
\end{equation}
Moreover, with $D:=M+C$,
\begin{equation}\label{mom:filtered-bounds}
 |u|,\ |u_x|,\ |A|,\ |B|,\ |\mathcal A|,\ |\mathcal B|\le D,
\end{equation}
and
\begin{equation}\label{mom:filtered-velocity-signs}
 \mathcal A(t,x)>0\Longrightarrow u(t,x)>0,
 \qquad
 \mathcal B(t,x)>0\Longrightarrow u(t,x)<0.
\end{equation}
\end{lemma}

\begin{proof}
Theorem~\ref{thm:entry} and \eqref{eq:J-label} give
\[
T^*=\infty\quad\Longrightarrow\quad J(t)=0
\quad\Longrightarrow\quad u_x(t,y)\ge-|u(t,y)|.
\]
At fixed $t$, we have
\[
\begin{aligned}
u(t,y)>0
&\ \Longrightarrow\
\partial_y(e^yu(t,y))=e^y(u_x+u)(t,y)\ge0,\\
u(t,y)<0
&\ \Longrightarrow\
\partial_y(e^{-y}u(t,y))=e^{-y}(u_x-u)(t,y)\ge0.
\end{aligned}
\]
Integrating on the corresponding sign interval excludes a first zero
to the right of a positive point or to the left of a negative point.
Thus
\begin{equation}\label{mom:spatial-comparison}
\begin{aligned}
u(t,x)>0&\ \Longrightarrow\
u(t,x+r)\ge u(t,x)e^{-r}>0\qquad(r\ge0),\\
u(t,x)<0&\ \Longrightarrow\
u(t,x-r)\le u(t,x)e^{-r}<0\qquad(r\ge0).
\end{aligned}
\end{equation}
For $u(t,x)\ne0$, use the sign $+$ when $u(t,x)>0$ and $-$ when
$u(t,x)<0$. By \eqref{mom:spatial-comparison} and \eqref{eq:energy-L2},
\[
\begin{aligned}
\frac12|u(t,x)|^2
&=|u(t,x)|^2\int_0^\infty e^{-2r}\dd r\\
&\le\int_0^\infty |u(t,x\pm r)|^2\dd r
\le\|u(t)\|_{L^2}^2\le4E_0.
\end{aligned}
\]
The bound also holds when $u(t,x)=0$. Hence
\[
|u(t,x)|\le M,\qquad u_x(t,x)\ge-|u(t,x)|\ge-M.
\]
For $W$ in \eqref{eq:UW}, the flow normalization \eqref{eq:flow}
and the estimates \eqref{sp:characteristic-identities} and
\eqref{eq:pressure-bound} give
\[
 \partial_tW(t,a)\le-W(t,a)^2+\frac32M^2,
 \qquad W(0,a)=u_0'(a).
\]
Comparison with the constant supersolution
$\max\{\sqrt{3/2}\,M,u_0'(a)\}$ gives
\[
 W(t,a)\le\max\{\sqrt{3/2}\,M,u_0'(a)\}\le C.
\]
Since $q(t,\cdot)$ is onto by Lemma~\ref{lem:characteristic} and
$C\ge M$, \eqref{eq:AB} and \eqref{eq:filtered-norms} give
\[
 |u|\le M,\qquad |u_x|\le C,\qquad
 |A|,|B|\le M+C=D,\qquad
 |\mathcal A|,|\mathcal B|\le D/2.
\]
By \eqref{mom:filtered-definitions}, \eqref{eq:AB}, and $u_x\ge-|u|$,
the positive convolution values yield the following implications, with
the last inequality in each coming from \eqref{mom:spatial-comparison}.
\[
\begin{aligned}
\mathcal A(t,x)>0
&\ \Longrightarrow\ \exists\,y_-<x:\quad
0<A(t,y_-)\le u(t,y_-)+|u(t,y_-)|\\
&\ \Longrightarrow\ u(t,y_-)>0,\qquad
u(t,x)\ge e^{y_--x}u(t,y_-)>0.
\end{aligned}
\]
\[
\begin{aligned}
\mathcal B(t,x)>0
&\ \Longrightarrow\ \exists\,y_+>x:\quad
0<B(t,y_+)\le |u(t,y_+)|-u(t,y_+)\\
&\ \Longrightarrow\ u(t,y_+)<0,\qquad
u(t,x)\le e^{x-y_+}u(t,y_+)<0.
\end{aligned}
\]
\end{proof}

\begin{remark*}
The implication $\mathcal A(t,x)>0\Rightarrow u(t,x)>0$ in
\eqref{mom:filtered-velocity-signs} uses the slope constraint.
For a nonzero, nonnegative $u\in C_c^\infty(\R)$ with
$\operatorname{supp}u\subset(-\infty,0)$, the spatial identity
\eqref{mom:filtered-algebra} gives
\[
u(0)=0,\qquad
\mathcal A(0)=3K_{2,-}u(0)-u(0)
=3\int_{-\infty}^0e^{2y}u(y)\dd y>0.
\]
\end{remark*}

\subsection{Local sign analysis}

\begin{lemma}[Positivity on an interval of positive momentum]\label{mom:positive-interval}
Suppose that the solution exists globally and $m_0>0$ on a nonempty
open interval $I$. Then there exist $t\ge0$ and $a\in I$ such that
$\mathcal A(t,q(t,a))>0$.
\end{lemma}

\begin{proof}
The proof is divided into four steps.
Suppose, to the contrary, that
\begin{equation}\label{mom:contradiction-assumption}
 \mathcal A(t,q(t,a))\le0\qquad(t\ge0,\ a\in I).
\end{equation}
By \eqref{eq:momentum-transport} and \eqref{eq:qa},
\[
m(t,q(t,a))=\frac{m_0(a)}{q_a(t,a)^3}>0
\qquad(t\ge0,\ a\in I).
\]
Since $u_0\not\equiv0$, \eqref{eq:energy-L2} gives $E_0>0$, and
\eqref{mom:global-estimates} gives $M>0$. Thus $D>0$ in
\eqref{mom:filtered-bounds}.
For fixed labels $a<b$ with $[a,b]\Subset I$, denote the image
endpoints by $Q_a(t),Q_b(t)$ and their distance by $L_{ab}(t)$.
We denote the cube of the local integral of $m_0^{1/3}$ by $\kappa_{ab}$.
\begin{equation}\label{mom:interval-notation}
 Q_a(t)=q(t,a),\qquad L_{ab}(t)=Q_b(t)-Q_a(t),\qquad
 \kappa_{ab}=\left(\int_a^bm_0(r)^{1/3}\dd r\right)^3>0.
\end{equation}
The change of variables in \eqref{eq:local-momentum-invariant}
follows from \eqref{eq:momentum-transport}.
\begin{equation}\label{mom:local-invariant}
 \begin{aligned}
 \int_{Q_a(t)}^{Q_b(t)}m(t,x)^{1/3}\dd x
 &=\int_a^b m(t,q(t,r))^{1/3}q_a(t,r)\dd r\\
 &=\int_a^b m_0(r)^{1/3}\dd r
 =\kappa_{ab}^{1/3}.
 \end{aligned}
\end{equation}
The cube in the definition of $\kappa_{ab}$ gives a lower bound for
the momentum integral via H\"older's inequality.

\proofstep{Step 1. Uniform bounds for the lengths of the image intervals.}
Using \eqref{mom:local-invariant}, H\"older's inequality,
$m=u-u_{xx}$ from \eqref{eq:momentum}, and \eqref{mom:filtered-bounds},
we obtain
\begin{equation}\label{mom:length-lower-calculation}
 \begin{aligned}
 \kappa_{ab}
 &=\left(\int_{Q_a}^{Q_b}m(t,x)^{1/3}\dd x\right)^3
 \le L_{ab}^2\int_{Q_a}^{Q_b}m(t,x)\dd x\\
 &=L_{ab}^2\left[\int_{Q_a}^{Q_b}u(t,x)\dd x
       +u_x(t,Q_a)-u_x(t,Q_b)\right]\\
 &\le D L_{ab}^2(L_{ab}+2).
 \end{aligned}
\end{equation}
If $L_{ab}(t)\le1$, this implies $\kappa_{ab}\le3D L_{ab}(t)^2$.
Including the case $L_{ab}(t)>1$, we obtain
\begin{equation}\label{mom:length-lower}
 L_{ab}(t)\ge\ell_{ab}:=
 \min\left\{1,\sqrt{\frac{\kappa_{ab}}{3D}}\right\}>0.
\end{equation}

Fix $a<b<c$ in $I$ and $t\ge0$. The convolution $\mathcal A=K_{2,-}A$ in
\eqref{mom:filtered-definitions} gives
\[
 e^{2(x-Q_a)}\mathcal A(t,x)
 =\int_{-\infty}^x e^{2(y-Q_a)}A(t,y)\dd y .
\]
Writing the weight as $z^2$ gives
\begin{equation}\label{mom:exponential-concavity}
 z=e^{x-Q_a},\qquad
 G(z)=-z^2\mathcal A(t,Q_a+\log z).
\end{equation}
For $1\le z\le z_c$, where $z_j=e^{L_{aj}}$ for $j=b,c$,
we have $G(z)\ge0$ by \eqref{mom:contradiction-assumption}.
At fixed $t$, put $x=Q_a+\log z$, so that
$\mathrm dx/\mathrm dz=1/z$. The product and chain rules, together with
\eqref{mom:filtered-algebra} and \eqref{eq:momentum-spatial}, give
\begin{equation}\label{mom:G-identities}
 \begin{aligned}
 G'(z)&=-2z\mathcal A(t,x)
        -z^2\mathcal A_x(t,x)\frac1z\\
      &=-z(2\mathcal A+\mathcal A_x)(t,x)
        =-zA(t,Q_a+\log z),\\
 G''(z)&=-A(t,x)-zA_x(t,x)\frac1z\\
       &=-(A+A_x)(t,x)=-m(t,Q_a+\log z)<0.
 \end{aligned}
\end{equation}
By \eqref{mom:G-identities} and \eqref{mom:filtered-bounds},
\[
0\le G(1)=-\mathcal A(t,Q_a)\le D,\qquad
G'(1)=-A(t,Q_a)\le D.
\]
By concavity, $G(z_b)\le G(1)+(z_b-1)G'(1)$ and $G'(z_b)\le G'(1)$.
Using \eqref{mom:local-invariant} and H\"older's inequality with
exponents $3$ and $3/2$, we obtain the following estimate. In the
product, we integrate the first factor by parts using
\eqref{mom:G-identities} and evaluate the second with the
antiderivative $-2z_c^{-1}\sqrt{(z_c-z)/z}$.
\begin{equation}\label{mom:weighted-mass-bound}
 \begin{aligned}
 \kappa_{bc}
 &=\left(\int_{z_b}^{z_c}
       m(t,Q_a+\log z)^{1/3}\frac{\dd z}{z}\right)^3\\
 &\le\left(\int_{z_b}^{z_c}(z_c-z)m(t,Q_a+\log z)\dd z\right)
       \left(\int_{z_b}^{z_c}
       \frac{\dd z}{z^{3/2}\sqrt{z_c-z}}\right)^2\\
 &=\frac{4(z_c-z_b)}{z_c^2z_b}
       \bigl[G(z_b)+(z_c-z_b)G'(z_b)-G(z_c)\bigr]\\
 &\le\frac{4(z_c-z_b)}{z_c^2z_b}\bigl[G(1)+(z_c-1)G'(1)\bigr]\\
 &\le4D\left(\frac1{z_b}-\frac1{z_c}\right)
 \le4De^{-L_{ab}(t)}.
 \end{aligned}
\end{equation}
Taking logarithms gives
\begin{equation}\label{mom:length-upper}
 L_{ab}(t)\le U_{ab}:=\log\frac{4D}{\kappa_{bc}}.
\end{equation}
Thus $U_{ab}$ is independent of $t$ once the extra label $c>b$ has
been fixed. For an outer interval $[a,d]\Subset I$, fix $e>d$ in $I$.
Applying \eqref{mom:length-upper} to $a<d<e$ gives
\[
 U_{ad}=\log\frac{4D}{\kappa_{de}},\qquad
 0<\ell_{ad}\le L_{ad}(t)\le U_{ad}<\infty\qquad(t\ge0).
\]
The same $U_{ad}$ bounds the image lengths of all subintervals of
$[a,d]$.

\proofstep{Step 2. A uniform strictly negative upper bound for \(\mathcal A\).}
Fix $a<b<c<d$ in $I$ and use the bound $U_{ad}$ from Step~1.
Keep the function $G$ from \eqref{mom:exponential-concavity} on
$[1,z_d]$, where $z_j=e^{Q_j-Q_a}$ for $j=b,c,d$.
By \eqref{mom:G-identities}, the integral Taylor formulas for
$G(1)$ and $G(z_d)$ about $z$ give
\[
\begin{aligned}
G(1)
&=G(z)-(z-1)G'(z)-\int_1^z(r-1)m(t,Q_a+\log r)\dd r,\\
G(z_d)
&=G(z)+(z_d-z)G'(z)
  -\int_z^{z_d}(z_d-r)m(t,Q_a+\log r)\dd r.
\end{aligned}
\]
Multiply the first equality by $z_d-z$ and the second by $z-1$.
Adding cancels $G'(z)$. Since $G(1),G(z_d)\ge0$, the resulting
estimate is
\begin{equation}\label{mom:concavity-interpolation}
 \begin{aligned}
 (z_d-1)G(z)
 &\ge (z_d-z)\int_{z_b}^{z}(r-1)m(t,Q_a+\log r)\dd r\\
 &\quad+(z-1)\int_z^{z_c}(z_d-r)m(t,Q_a+\log r)\dd r\\
 &\ge (z_b-1)(z_d-z_c)\int_{z_b}^{z_c}m(t,Q_a+\log r)\dd r
 \qquad(z_b\le z\le z_c).
 \end{aligned}
\end{equation}
The length bounds \eqref{mom:length-lower} and
\eqref{mom:length-upper}, together with $e^r-1\ge r$ for $r\ge0$,
give
\[
 z_d\le e^{U_{ad}},\qquad z_b-1\ge\ell_{ab},\qquad
 z_d-z_c=z_c(e^{Q_d-Q_c}-1)\ge\ell_{cd}.
\]
By \eqref{mom:length-lower-calculation},
\[
 \begin{aligned}
 \int_{z_b}^{z_c}m(t,Q_a+\log r)\dd r
 &=\int_{Q_b}^{Q_c}e^{x-Q_a}m(t,x)\dd x\\
 &\ge\int_{Q_b}^{Q_c}m(t,x)\dd x
 \ge\frac{\kappa_{bc}}{L_{bc}(t)^2}
 \ge\frac{\kappa_{bc}}{U_{ad}^2}.
 \end{aligned}
\]
Thus \eqref{mom:concavity-interpolation} yields
\begin{equation}\label{mom:strict-negative}
 \begin{gathered}
 -\mathcal A(t,x)=\frac{G(z)}{z^2}\ge e^{-2U_{ad}}G(z)\ge\delta
 \qquad(Q_b\le x\le Q_c),\\
 \delta:=\frac{e^{-2U_{ad}}\ell_{ab}\ell_{cd}\kappa_{bc}}
 {(e^{U_{ad}}-1)U_{ad}^2}>0.
 \end{gathered}
\end{equation}
Any fixed compact subinterval of $I$ can serve as $[b,c]$ by choosing
outer labels $a,d$ inside $I$.

\proofstep{Step 3. Excluding a positive value of \(\mathcal B\).}
Suppose that $\mathcal B(t_0,q(t_0,a_0))>0$ for some $a_0\in I$.
By continuity, choose an open interval $I'\Subset I$ containing $a_0$
on which $\mathcal B(t_0,q(t_0,a))>0$. By
\eqref{mom:weighted-monotonicity},
\[
e^{-q(t,a)}\mathcal B(t,q(t,a))
\ge e^{-q(t_0,a)}\mathcal B(t_0,q(t_0,a))>0
\qquad(t\ge t_0,\ a\in I').
\]
By \eqref{mom:filtered-velocity-signs}, \eqref{eq:momentum-transport},
and \eqref{eq:momentum},
\[
\begin{gathered}
u(t,q(t,a))<0<m(t,q(t,a)),\\
u_{xx}(t,q(t,a))=(u-m)(t,q(t,a))<0
\qquad(t\ge t_0,\ a\in I').
\end{gathered}
\]
Fix $a<b<c<d$ in $I'$.
By \eqref{mom:length-lower} and \eqref{mom:length-upper},
\begin{equation}\label{mom:ratio-lower}
 \frac{L_{cd}(t)}{L_{ab}(t)}\ge\frac{\ell_{cd}}{U_{ad}}>0
 \qquad(t\ge t_0).
\end{equation}
For $x\in[Q_a,Q_b]$ and $y\in[Q_c,Q_d]$, the sign of $u$
and \eqref{mom:length-lower-calculation} give
\[
\begin{aligned}
u_x(t,y)-u_x(t,x)
&=\int_x^y(u-m)(t,r)\dd r
 \le-\int_{Q_b}^{Q_c}m(t,r)\dd r\\
&\le-\frac{\kappa_{bc}}{L_{bc}(t)^2}
 \le-\frac{\kappa_{bc}}{U_{ad}^2}
 \qquad(t\ge t_0).
\end{aligned}
\]
By \eqref{eq:interval-log-derivative}, the logarithmic derivatives
of $L_{ab}$ and $L_{cd}$ are the spatial averages of $u_x$ on the
respective intervals. Thus
\begin{equation}\label{mom:ratio-decay}
 \begin{aligned}
 \frac{\mathrm d}{\mathrm dt}
       \log\frac{L_{cd}(t)}{L_{ab}(t)}
 &=\frac{1}{L_{ab}(t)L_{cd}(t)}
 \int_{Q_a}^{Q_b}\int_{Q_c}^{Q_d}
   [u_x(t,y)-u_x(t,x)]\dd y\dd x\\
 &\le-\frac{\kappa_{bc}}{U_{ad}^2}<0 .
 \end{aligned}
\end{equation}
Integrating from $t_0$ gives
\[
\frac{L_{cd}(t)}{L_{ab}(t)}
\le\frac{L_{cd}(t_0)}{L_{ab}(t_0)}
   \exp\!\left(-\frac{\kappa_{bc}}{U_{ad}^2}(t-t_0)\right)
\longrightarrow0,
\]
contrary to \eqref{mom:ratio-lower}. Hence
\begin{equation}\label{mom:B-nonpositive}
 \mathcal B(t,q(t,a))\le0\qquad(t\ge0,\ a\in I).
\end{equation}

\proofstep{Step 4. A contradiction from the decay of \(\mathcal B\).}
Fix $a<b<c<d$ in $I$.
Apply \eqref{mom:strict-negative} to $[a,d]$ as the inner label
interval, with a fixed label on each side. By
\eqref{mom:filtered-bounds}, we have
\[
0<\delta\le-\mathcal A(t,q(t,r))\le D
\qquad(t\ge0,\ r\in[a,d]).
\]
The first identity in \eqref{mom:weighted-monotonicity} gives
\[
\begin{aligned}
\frac{\mathrm d}{\mathrm dt}
 \log\!\left(-e^{q(t,r)}\mathcal A(t,q(t,r))\right)
&=\frac12\mathcal A(t,q(t,r))
 +\frac32\frac{K_-(\mathcal A^2)(t,q(t,r))}{\mathcal A(t,q(t,r))}\\
&\le-\frac{\delta}{2}.
\end{aligned}
\]
Using $q(0,r)=r$ from \eqref{eq:flow} and integrating gives
\begin{equation}\label{mom:left-drift}
 q(t,r)-r
 \le\log\frac{-\mathcal A(0,r)}{-\mathcal A(t,q(t,r))}
 -\frac\delta2t
 \le\log\frac D\delta-\frac\delta2t
 \qquad(r\in[a,d]).
\end{equation}
By \eqref{mom:B-nonpositive}, \eqref{mom:weighted-monotonicity},
and \eqref{mom:left-drift},
\begin{equation}\label{mom:B-decay}
 \begin{aligned}
 0\le-\mathcal B(t,q(t,r))
 &\le[-\mathcal B(0,r)]e^{q(t,r)-r}\\
 &\le\frac{D^2}{\delta}e^{-\delta t/2}
 \qquad(r\in[a,d]).
 \end{aligned}
\end{equation}

For the momentum integral, the convolution $\mathcal B=K_{2,+}B$
in \eqref{mom:filtered-definitions} gives
\[
 e^{2(Q_d-x)}\mathcal B(t,x)
 =\int_x^\infty e^{2(Q_d-y)}B(t,y)\dd y .
\]
Let $\bar{z}:=e^{Q_d-x}$ and set $H(\bar{z}):=-\bar{z}^2\mathcal B(t,Q_d-\log \bar{z})$.
The endpoints are $\bar{z}_j=e^{Q_d-Q_j}$, with
$1= \bar{z}_d<\bar{z}_c<\bar{z}_b<\bar{z}_a\le e^{U_{ad}}$.
At fixed $t$, put $x=Q_d-\log\bar z$, so that
$\mathrm dx/\mathrm d\bar z=-1/\bar z$. The product and chain rules,
together with \eqref{mom:filtered-algebra} and
\eqref{eq:momentum-spatial}, give
\begin{equation}\label{mom:H-identities}
\begin{aligned}
 H'(\bar z)&=-2\bar z\mathcal B(t,x)
 -\bar z^2\mathcal B_x(t,x)\left(-\frac1{\bar z}\right)\\
 &=\bar z(\mathcal B_x-2\mathcal B)(t,x)
   =-\bar zB(t,Q_d-\log\bar z),\\
 H''(\bar z)&=-B(t,x)
 -\bar zB_x(t,x)\left(-\frac1{\bar z}\right)\\
 &=(B_x-B)(t,x)=m(t,Q_d-\log\bar z)>0.
\end{aligned}
\end{equation}
Thus $H$ is convex, and \eqref{mom:B-decay} gives
\[
 0\le H(\bar{z})\le\frac{D^2e^{2U_{ad}}}{\delta}e^{-\delta t/2}
 \qquad(1\le \bar{z}\le \bar{z}_a).
\]
By \eqref{mom:length-lower}, the two outer gaps satisfy
\[
 \bar{z}_a-\bar{z}_b=\bar{z}_b(e^{Q_b-Q_a}-1)\ge\ell_{ab},\qquad
 \bar{z}_c-1=e^{Q_d-Q_c}-1\ge\ell_{cd}.
\]
By \eqref{mom:length-lower-calculation} and \eqref{mom:length-upper},
\begin{equation}\label{mom:mass-lower}
 \int_{Q_b}^{Q_c}m(t,x)\dd x
 \ge\frac{\kappa_{bc}}{L_{bc}(t)^2}
 \ge\frac{\kappa_{bc}}{U_{ad}^2}>0.
\end{equation}
On the other hand, using \eqref{mom:H-identities}, changing variables,
and comparing the derivatives of $H$ with the outer secant slopes give
\begin{equation}\label{mom:mass-upper}
 \begin{aligned}
 \int_{Q_b}^{Q_c}m(t,x)\dd x
 &=\int_{\bar{z}_c}^{\bar{z}_b}H''(\bar{z})\frac{\dd \bar{z}}{\bar{z}}
 \le H'(\bar{z}_b)-H'(\bar{z}_c)\\
 &\le\frac{H(\bar{z}_a)-H(\bar{z}_b)}{\bar{z}_a-\bar{z}_b}
    -\frac{H( \bar{z}_c)-H(1)}{\bar{z}_c-1}\\
 &\le\left(\ell_{ab}^{-1}+\ell_{cd}^{-1}\right)
       \frac{D^2e^{2U_{ad}}}{\delta}e^{-\delta t/2}\longrightarrow0.
 \end{aligned}
\end{equation}
This contradicts \eqref{mom:mass-lower} and proves the lemma.
\end{proof}

\subsection{Characterization by momentum signs}

\begin{proof}[Proof of Theorem~\ref{thm:momentum}]
First suppose that there exist $a_0<b_0$ with
$m_0(a_0)>0>m_0(b_0)$, and assume, to the contrary, that the solution
exists globally. By Lemma~\ref{lem:momentum-identities}, $m_0$ is
continuous. Choose bounded open intervals $I_+$
and $I_-$, with $I_+$ entirely to the left of $I_-$, on which
$m_0>0$ and $m_0<0$, respectively.
Lemma~\ref{mom:positive-interval} gives $\alpha\in I_+$ and
$t_A\ge0$ such that
\[
 \mathcal A(t_A,q(t_A,\alpha))>0.
\]

By \eqref{eq:reflection}, the reflected solution
$\widetilde u(t,x)=-u(t,-x)$ has positive initial momentum on $-I_-$,
and $\widetilde{\mathcal A}(t,x)=\mathcal B(t,-x)$.
Applying the same lemma to this solution gives $\beta\in I_-$ and
$t_B\ge0$. Since $\widetilde q(t,-\beta)=-q(t,\beta)$,
\[
 \mathcal B(t_B,q(t_B,\beta))>0.
\]
Let $T=\max\{t_A,t_B\}$. By \eqref{mom:weighted-monotonicity},
both quantities remain strictly positive at time $T$.
Hence \eqref{mom:filtered-velocity-signs} gives
\[
 u(T,q(T,\alpha))>0>u(T,q(T,\beta)).
\]
Since $\alpha<\beta$, Lemma~\ref{lem:characteristic} gives
$q(T,\alpha)<q(T,\beta)$. Applying \eqref{mom:spatial-comparison}
with $t=T$, $x=q(T,\alpha)$, and
$r=q(T,\beta)-q(T,\alpha)>0$, we obtain
\[
 u(T,q(T,\beta))
 \ge e^{-[q(T,\beta)-q(T,\alpha)]}u(T,q(T,\alpha))>0.
\]
This contradicts $u(T,q(T,\beta))<0$, so the maximal lifespan is finite.

Conversely, suppose no such pair exists. If both signs occur, every
negative point lies to the left of every positive point. The number
$x_*=\sup\{x\in\R:m_0(x)<0\}$ is therefore finite, and continuity at
$x_*$ gives \eqref{eq:global-sign}.
The global-existence results of Liu and Yin
\cite[Theorems~3.1 and~3.2]{LiuYin2006} give $T^*=\infty$ in
the one-sign and the ordered-sign cases.
Finally, the wave-breaking assertion \eqref{eq:breaking} follows from
Theorem~\ref{thm:entry}.
\end{proof}

\medskip
\noindent\textbf{Acknowledgments.}
Zhaoyang Yin was partially supported by the National Natural Science Foundation of China (Grant No.~12571261).
Weikui Ye was supported by the National Natural Science Foundation of China (Grant No.~12401277).

\medskip
\noindent\textbf{Conflict of interest.}
The authors declare that they have no conflict of interest.

\medskip
\noindent\textbf{Data availability statement.}
No datasets were generated or analysed during the current study.

\end{document}